\documentclass[a4paper,reqno, 10pt]{amsart}

\usepackage{amsthm}
\usepackage{amsmath}
\usepackage{amssymb}
\usepackage[dvipsnames]{xcolor}
    \definecolor{awesome}{rgb}{1.0, 0.13, 0.32}
    \definecolor{Mulberry}{rgb}{0.66, 0.24, 0.58}
\usepackage{graphicx}
\usepackage{tikz}
\usetikzlibrary{cd}
\usetikzlibrary{decorations.pathreplacing,calligraphy,calc}
\usepackage{fullpage}
\usepackage[english]{babel}
\usepackage{bm}
\usepackage{txfonts}
\usepackage{comment}
\usepackage{enumitem}
\usepackage{dsfont}
    
\usepackage{microtype}
\usepackage[colorlinks = true, linkbordercolor = {white}, urlcolor={purple}, linkcolor={purple}, citecolor = {purple}]{hyperref}
\usepackage{cleveref}
\usepackage[foot]{amsaddr}
\usepackage{marginnote}
\usepackage{mathrsfs}
\usepackage{subcaption}
\usepackage{array}
\usepackage{multirow}
\usepackage{transparent}
\usepackage{float}
\usepackage{wrapfig}
\usepackage{mathtools}

\numberwithin{equation}{section}

\renewcommand{\emph}{\textsl}

\renewcommand{\textit}{\textsl}

\title{Topological notions of Shadowing under semi-conjugacy }

\author[A.\,Bura]{A.\,Bura}
\author[C.\,Good]{C.\,Good}

\address[A.\,Bura]{School of Mathematics, University of Birmingham, UK}
\address[C.\,Good]{School of Mathematics, University of Birmingham, UK }
    
    \theoremstyle{definition}
    \newtheorem{definition}{Definition}[section]

    \newtheorem{theorem}{Theorem}[section]

\theoremstyle{definition}
    
\theoremstyle{remark}

\begin{document}

\begin{abstract}
    We study various notions of shadowing and investigate whether they are preserved under semi-conjugacy and inverse limits. In particular we study Inverse shadowing, $H$-shadowing, Periodic shadowing, Ramsey shadowing, Density shadowing and Ergodic shadowing.
\end{abstract}

\maketitle

\section{Introduction}

Let $(X,d)$ be a compact metric space and  $f : X \to X$ be a continuous function. Given a point $z \in X$, we define the orbit of $z$ as the set $\{z, f(z),f^2(z), \dots \}$. Let $\delta >0 $. A $\delta$\emph{-pseudo-orbit} is a sequence $(x_i)_{i\in \mathbb{N}_0}$ such that $\mathrm{d}(f(x_{i}),x_{i+1}) <\delta$ for all $i \in \mathbb{N}_{0}$. Intuitively a $\delta$-pseudo orbit is a sequence that behaves like an orbit, however at every stage of the iteration there is an error of at most $\delta$.   Given $\epsilon >0$, a sequence $(y_i)_{i \in \mathbb{N}_{0}}$ in $X$ is said to \emph{$\epsilon$-shadow} the sequence $(x_i)_{i\in \mathbb{N}_{0}}$ if $\mathrm{d}(y_i,x_i)<\epsilon$ for all $i \in \mathbb{N}_{0}$. We say that $f$ has \emph{shadowing} if every $\epsilon >0 $, there exists a $\delta > 0$, such that every $\delta$-pseudo-orbit is $\epsilon$-shadowed by a true orbit. 

Shadowing has become an important tool in the  study of dynamical systems. Anosov's shadowing lemma shows that diffeomorphisms restricted to a hyperbolic set exhibit shadowing \cite{anosov1967geodesic}. This result has been used to study the stability of hyperbolic dynamical systems. Bowen used shadowing to study Axiom A diffeomorphisms and their $\omega$-limit sets \cite{BOWEN1975333}.

 Shadowing has also emerged as an important tool in applied contexts, particularly in numerical modelling. Computers operate with finite precision, hence the iteration of a function typically introduces small rounding errors. The resulting numerical trajectory can therefore be viewed as a pseudo-orbit rather than an exact orbit. If the dynamical system possesses the shadowing property, then this pseudo-orbit remains close to a genuine orbit. Thus, despite the presence of numerical errors, the computed trajectory can still be interpreted as representing a valid behaviour of the underlying dynamical system. For this reason, shadowing has been studied extensively in the context of numerical modelling and the reliability of numerical simulations \cite{HAYES2005299}.

In \cite{Good} it was shown that shadowing can be formulated in purely topological terms and therefore does not require the underlying space to be metrizable. Building on this observation, they introduced the almost lifting property, which characterizes when a factor map preserves the shadowing property.
 
A factor map between two dynamical systems $(X,f),(Y,g)$ is a continuous surjection, $\psi: X \longrightarrow Y$, such that the following holds,  $\psi \circ f$ = $g \circ \psi$.  The authors also showed that shadowing was preserved under inverse limits provided it satisfied the Mittag-Leffler condition.

 Since shadowing was first introduced, various other notions of shadowing have been developed. These include inverse shadowing, orbital shadowing and ergodic shadowing. Work has been done on studying whether these differing notions of shadowing are preserved under various actions such as factor maps and inverse limits \cite{MT}.

 In this paper, we expand upon the work of \cite{MT} by considering several other notions of shadowing and addressing some of the unresolved questions in the paper. In particular, we  look at inverse shadowing, h-shadowing, periodic shadowing, Ramsey shadowing, density shadowing and ergodic shadowing. We characterize the property that determines whether they are preserved under factor maps and investigate whether they are preserved under inverse limits, and the hyperspace topology.

\section{Preliminaries}

\begin{definition}
    A dynamical system is a compact Hausdorff space $X$, with a continuous function $f : X \to X$.
\end{definition}

The orbit of a point $x \in X$ is defined to be the set $\{x, f(x),f^2(x),\dots \}$.
Given two dynamical systems $(X,f),(Y,g)$, a factor map from $(X,f)$ to $(Y,g)$  is a continuous surjection, $\psi$ , such that the following holds,  $\psi \circ f$ = $g \circ \psi$. $\psi$ is often also called a semi-conjugacy and we say that $(X,f) \text{ and }(Y,g)$ are semi-conjugate. If $\psi$ is also a homeomorphism, we call it a conjugacy, and say that  $(X,f) \text{ and }(Y,g)$ are conjugate.

\begin{definition}
    
Given a dynamical system $(X,f)$, let $2^X$ denote the set of non-empty compact subsets of $X$. Given a subset $B$ of $X$ define $B_{\delta} \coloneqq  \{z: z \in B(x,\delta), x \in B  \}$. For subsets $A,B$ of $X$, we define the Hausdorff metric by $d_H(A,B) \coloneqq \inf\{\delta : A \subseteq B_\delta \text{ and }  B\subseteq A_\delta \} $. . This forms a compact topological space under the Hausdorff metric. Then $f$ induces a continuous function, $F$ from $2^X$ to $2^X$ defined by sending $B \in 2^X $ to $f(B)$. Hence $(2^X,F)$ is also a dynamical system.
\end{definition}

\bigskip

\begin{definition}

Consider a collection of dynamical systems $(X_n,F_n)_{n \in \mathbb{N}_0}$. Suppose for all $m \leq n$ we have maps $\psi^{n}_{m}: X_n \longrightarrow X_m  $ which satisfy: 

\begin{itemize}
    \item $ \psi^{n}_{n} $ is the identity.

    \item  For all $ l \leq m \leq n$ we have that $\psi^{n}_{l}= 
           \psi^{m}_{l}   \circ \psi^{n}_{m} $.

    \item{$F_m \circ  \psi^{n}_{m}$ = $ \psi^{n}_{m}   \circ F_n$ }

    \end{itemize}

We call this an inverse system $(\psi^{n}_{m},X_n)$ and its inverse limit space is given by :

$   \varprojlim\{ \psi^{n}_{m}, X_n  \}= \{ (x_i) \in \prod X_n \ : \forall m \leq n \   \psi^{n}_{m}(x_n)=x_m  \}$

The inverse limit is the dynamical system given by $ ( \varprojlim\{ \psi^{n}_{m}, X_n  \},F^*   )$ where $F_{}^*(\langle x_n \rangle )= (F_{1}(x_1),F_{2} \\ (x_2),F_{3}(x_3), \dots   ) $

\end{definition}

An inverse system of dynamical systems satisfies the Mittag-Leffler condition provided that for all $l$, we can find an $m \geq l$ such that for all $n \geq m$ we have that $\psi^{n}_{l}(X_n)= \psi^{m}_{l}(X_m) $. Clearly if each $\psi$ is surjective, the Mittag-Leffler condition is satisfied.

\bigskip

\begin{definition}[Shadowing]
Let $\delta >0 $. A sequence $(x_i)_{i\in \mathbb{N}_0}$  is called a $\delta$\emph{-pseudo-orbit} if $\mathrm{d}(f(x_{i}),x_{i+1}) <\delta$ for all $i \in \mathbb{N}_0$. Given $\epsilon >0$, a sequence $(y_i)_{i \in \mathbb{N}_{0}}$ in $X$ is said to \emph{$\epsilon$-shadow} the sequence $(x_i)_{i\in \mathbb{N}_{0}}$ if $\mathrm{d}(y_i,x_i)<\epsilon$ for all $i \in \mathbb{N}_{0}$. We say that $f$ has \emph{shadowing} if for every $\epsilon >0 $, there exists a $\delta > 0$, such that if $(x_i)_{i\in \mathbb{N}_0}$ is a $\delta$-pseudo-orbit, then there exists a $z$ such that $(x_i)_{i\in \mathbb{N}_0}$ is $\epsilon$-shadowed by a true orbit.

\end{definition}

The definition given for shadowing is given in a metric space. However in \cite{Good}, Good and Meddaugh showed that an equivalent definition of shadowing could be formulated on topological spaces.

\begin{definition}
Let $(X,f)$ be a dynamical system, where $X$ is a topological compact space. Let $FOC(X)$ denote the set of finite open covers of $X$. Let $U \in FOC(X)$. We say a sequence  $(x_i)_{i\in \mathbb{N}_0}$ is a $U$-pseudo orbit provided that for all $i \in \mathbb{N}_0$, there exists a $U_{i+1} \in U$ such that $f(x_i),x_{i+1} \in U_{i+1}$.

We say a sequence $(x_i)_{i\in \mathbb{N}_0}$ is $U$-shadowed by a sequence $(y_i)_{i\in \mathbb{N}_0}$ if for all  $i \in \mathbb{N}_0$, there exists a $U_i \in U$ such that $x_i,y_i \in U_i$.
We say $(X,f)$ has topological shadowing if for all $U \in FOC(X)$, there exists $Y \in FOC(X)$ such that if $(x_i)_{i\in \mathbb{N}_0}$ is a $U$-pseudo orbit , then there a $z \in X$ such that the orbit of $z$ $Y$-shadows $(x_i)_{i\in \mathbb{N}_0}$.

\end{definition}

\begin{theorem} [{\cite[Lemma 6]{Good}}]

Let $X$ be a compact metric space. A dynamical system $(X,f)$ has shadowing  if and only if for every $U \in FOC(X)$, there exists a $V \in FOC(X)$ such that every $V$-pseudo orbit is $U$-shadowed by a point $z \in X$. 
    
\end{theorem}

 Analogously to \cite{Good}, for each metric definition of an alternative notion of shadowing, we also provide the corresponding topological definition. We will show that these definitions are equivalent to their corresponding metric definition, provided the underlying metric space $X$ is compact.

\section{Inverse shadowing}

\begin{definition}\label{def:3} Inverse shadowing- Metric definition:

Let $(X,f)$ be a dynamical system, and let $X^\omega$ denote the set of all infinite sequences of $X$. For $\delta>0$, let $\Phi_f(\delta) \subseteq X^\omega$ denote the set of all sequences given by a $\delta$-pseudo orbit. A $\delta$-method for $f$ is a map $\rho: X \rightarrow \Phi_f(\delta)$, which sends every $x \in X$ to a pseudo orbit sequence which begins with $x$. We let $\mathcal{T}(f,\delta)$ denote the set of all $\delta$-methods of $f$.

A dynamical system  $(X,f)$ has inverse shadowing if for all $\epsilon>0$ there exist $\delta>0$ such that for all $x$$\in$$X$ and all $\delta_{f}$ $\in$ $\mathcal{T}(f,\delta)$, there exists $y$ such that, $\delta_{f}(y)$ $\epsilon$-shadows the orbit of $x$.

\end{definition}

Under the standard notion of shadowing, we are interested in whether pseudo-orbits can be approximated by a real orbit. Inverse shadowing consider the  converse problem: whether any orbit can be approximated using a suitable pseudo-orbit from any $\delta$-method of $f$. Inverse shadowing arises when studying the structural stability and hyperbolicity of various dynamical systems \cite{lee2009continuous}. In particular, it has been shown that on compact manifolds, hyperbolic homeomorphisms have inverse shadowing under certain conditions on the set of $\delta$-methods \cite{Pilyugin2002Inverse}.

\begin{definition}\label{def:4} Inverse shadowing- Topological definition:

A dynamical system  $(X,f)$ has inverse shadowing if for all $V_{X}$ $\in$ $FOC(X)$, there exists an $W_{X}$ $\in$ $FOC(X)$, such that  for all $x$$\in$$X$ and all $\delta_{f}$ $\in$ $\mathcal{T}(f,W_{X})$, there exists an 
$y$ such that, $\delta_{f}(y)$ $V_{X}$-shadows $x$.

\end{definition}

We use the following well-known theorem to prove the equivalence between the topological and metric notions of inverse-shadowing in a compact metric space. 
\begin{theorem}(Lebesgue's number lemma)
    Let $X$ be a compact metric space. Let $U$ be an open cover of $X$. There exists a $\delta>0$ such that for every $x \in X$, there exists a $U_k \in U$ such that $B(x,\delta) \subset U_k$.  We call this $\delta$ the Lebesgue number of $U$.
\end{theorem}

\begin{theorem}
    Let $(X,f)$ be a dynamical system where $X$ is a compact metric space. The topological and metric definitions of inverse shadowing are equivalent.
\end{theorem}

\begin{proof}
    Suppose the topological definition of inverse shadowing holds in $(X,f)$. Let $\epsilon>0$. Since $X$ is compact we can choose an open cover $U$ of $X$ which consists of a finite number of $\frac{\epsilon}{4}$ balls. Choose $V \in FOC(X)$ which satisfies \Cref{def:4} with respect to $U$. Let $\delta>0$ denote the Lebesgue number of $V$. Pick an arbitrary $x \in X$ and an arbitrary  $\delta_{f}$ $\in$ $\mathcal{T}(f,\delta)$. Then   $\delta_{f}$ is a map $\rho: X \rightarrow \Phi_f(\delta)$, which sends every $y \in X$ to a $\delta$-pseudo orbit sequence which begins with $y$. For all $i \in \mathbb{N}_0$ we have $d(f(y_i),y_{i-1})<\delta$. As $\delta$ is the Lebesgue number of $V$, we know that there exists $V_k \in V$ such that $f(y_i),y_{i-1} \in V_k$. Hence $\Phi(y)$ is also a $V$-pseudo orbit which starts with $y$. Therefore $\delta_f$ is also a member of $\mathcal{T}(f,V)$. Choose a $y$ which satisfies \Cref{def:4} with respect to $x$ and $\delta_f$. As we have topological inverse shadowing, we know that for all $i \in \mathbb{N}_0$, there exists a $U_j \in U$ such that $f^i(x)\text{ and }(\delta_f(y))_i$ both belong to $U_j$. Since $U$ consists of $\frac{\epsilon}{4}$-balls, we know that $d(f^i(x), (\delta_f(y))_i  )<\frac{\epsilon}{2}$. Hence the metric definition of inverse shadowing is satisfied.

    Now we consider the case where the metric definition holds. Consider an finite open cover $U$ of $X$. Let $\epsilon>0$ denote the Lebesgue number of $U$. Choose $\delta>0$ which satisfies \Cref{def:3} with respect to $\epsilon$. As $X$ is compact we can choose a finite open cover $V$ of $X$ which consists of $\frac{\delta}{2}$-balls. Pick an arbitrary $x \in X$ and an arbitrary  $\delta_{f}$ $\in$ $\mathcal{T}(f,V)$.
    
    Then   $\delta_{f}$ is a map $\rho: X \rightarrow \Phi_f(V)$, which sends every $y \in X$ to a $V$-pseudo orbit sequence which begins with $y$. For all $i \in \mathbb{N}_0$, there exists $V_k \in V$ such that $f(y_i),y_{i-1} \in V_K$. $V$ consists of $\frac{\delta}{2}$-balls, so we know that $d(f(y_i),y_{i-1}) < \delta$.  Hence $\Phi(y)$ is also a $\delta$-pseudo orbit which starts with $y$. Therefore $\delta_f$ is also a member of $\mathcal{T}(f,\delta)$. Choose $y$ which satisfies \Cref{def:3}. For $i \in \mathbb{N}_0$ we have $d(f^i(x),y_i)<\epsilon$ . Hence there exists a $U_j$ such that $f^i(z),x_i \in U_j $. Hence the topological definition of inverse shadowing is also satisfied.

\end{proof}

\begin{definition}Almost lifting delta method property (ALDM)- Metric definition:

A factor map $\psi:X \longrightarrow Y $ has almost lifting delta method property, if for all $\epsilon>0$ and $\eta>0$ there exist $\delta>0$, such that if $\delta_{g}$ $\in$ $\mathcal{T}(g,\delta)$, then there exists $\delta_{f}$ $\in$ $\mathcal{T}(f,\eta)$ such that $\psi(\delta_{f})$ $\epsilon$-shadows $\delta_{g}$.

When we say that $\psi(\delta_{f})$ $\epsilon$-shadows $\delta_{g}$, we mean that for all $x$ $\in$ $X$ we have $\psi(\delta_{f}(x))$ $\epsilon$-shadows some element of $\delta_{g}$.
    
\end{definition}

\begin{definition}Almost lifting delta method property (ALDM)- Topological definition:

$\psi:X \longrightarrow Y $ has almost lifting delta method property, if for all $V_{Y}$ $\in$ $FOC(Y)$ and $V_{X}$ $\in$ $FOC(X)$ there exist $W_{Y}$ $\in$ $FOC(Y)$, such that if $\delta_{g}$ $\in$ $\mathcal{T}(g,W_{Y})$, then there exists $\delta_{f}$ $\in$ $\mathcal{T}(f,V_{X})$ such that $\psi(\delta_{f})$ $V_{Y}$-shadows $\delta_{g}$.

\end{definition}

It is straightforward to show that the topological and metric definitions of ALDM are equivalent, hence we omit the proof.

\begin{theorem}
    Suppose $\psi$ is a factor map from  $(X,f)$ to $(Y,g)$. We will prove that if $(Y,g)$ has inverse shadowing, then $\psi$ must have ALDM.

\end{theorem}

\begin{proof}

Let $\epsilon > 0$ and $\eta>0$. Choose $\delta >0$ such that $(Y,g)$  witnesses  inverse shadowing with respect to $\epsilon$. Let $y$ $\in$ $Y$ be arbitrary. Let $\delta_{g}$ $\in$ $\mathcal{T}(g,\delta)$ be arbitrary. For each $x$ $\in$ $X$, consider its orbit. This clearly define a $\eta$-pseudo orbit for each $x$ $\in$ $X$. So consider $\delta_{f}$ $\in$ $\mathcal{T}(f,\eta)$ defined by the orbits on $X$. Let $x \in X$ be arbitrary.  Now consider $\delta_{f}(x)$. Then clearly $\psi(\delta_{f}(x))$ is the orbit of $\psi(x)$. But since $(Y,g)$ has inverse shadowing, we know there exists a $y$ such that $\delta_{g}(y)$ $\epsilon$-shadows the orbit of $\psi(x)$, and hence $\psi(\delta_{f}(x))$ $\epsilon$-shadows  $\delta_{g}(y)$. Hence $\psi$ has the ALDM property.

\end{proof}

\begin{theorem}

Let $(X,f)$ and $(Y,g)$ be two dynamical systems. suppose $(X,f)$ has inverse shadowing and $\psi:X \longrightarrow Y$ has the ALDM property, then $(Y,g)$ has inverse shadowing.
\end{theorem}

\begin{proof}

Fix an open cover $U_{Y}$ of $Y$. Let $W_{Y}$ be an open cover of $Y$ such that if $W$,$W'$ $\in$  $W_{Y}$ and $W$ $\cap$ $W'$$\neq$$\emptyset$,  there exist $U$$\in$$U_{Y}$ such that  $W$ $\cup$ $W'$ $\subset U $. Let $U_{X}$ be an open cover of X such that $\psi(U_{X})$ refines $W_Y$. Since $(X,f)$ has inverse shadowing , let $V_{X}$ witness this with respect to $U_{X}$. Since $\psi$ has ALDM, let $V_Y$ witness this with respect to $W_Y$ and $V_X$.

\bigskip

Let $y \in Y$. Let $\delta_{g} \in \mathcal{T}(g,V_{Y})$ be arbitrary. Consider $z \in \psi^{-1}(y)$. By inverse shadowing we know that for all $\delta_{f}$ $\in$ $\mathcal{T}(f,V_{X})$, there exists a $x' \in X$, such that $\delta_{f}(x')$ $U_{X}$-shadows $z$. Take $\delta_{f^*} \in \mathcal{T}(f,V_{X})$, such that it satisfies the conditions of the ALDM property. Hence, we know that z is $U_{X}$-shadowed by $\delta_{f^*}(x)$ for some $x \in X$. Since $\psi(U_{X})$ refines $W_Y$. We must have that for all $i$, $\psi(f^i(z))$ and $(\psi(\delta_{f^*}(x))_{i} $ are both in $W_{Y_{k}}$, for some $k$.  But notice that $\psi(f^i(z))=g^i(\psi(z))=g^i(y)$. Hence for all $i$, $g^i(y))$ and $(\psi(\delta_{f^*}(x))_{i} $ are both in $W_{Y_{k}}$.

\bigskip
But by the ALDM property $\psi(\delta_{f^{*}}(x))$ $W_{Y}$-shadows $\delta_{g}(y')$ for some $ y' \in Y$. Hence for all $i$, $(\psi(\delta_{f^{*}}(x)))_i$ and $(\delta_{g}(y'))_i$ belong to $W_{Y_{t}}$ for some $t$. So for all $i$ we must have that $W_{Y_{k}}$ and $W_{Y_{t}}$ both contain $(\psi(\delta_f^*(x)))_i$, and hence their intersection is non empty. But by construction of $W_{Y}$, we know that there exists a $s$ such that $g^i(y))$ and $(\delta_{g}(y'))_i$ are both contained in $U_{Y_{s}}$. Hence $y$ is $U_Y$-shadowed by  $\delta_{g}(y')$. And so we have shown that $(Y,g)$ has inverse shadowing.

\end{proof}

In \cite{MT} it was shown that that inverse shadowing is preserved under inverse limits and the induced hyperspace.

\section{$H$-shadowing}

\begin{definition}\label{def:$H$-metric}$H$-shadowing - Metric definition:
    We say that $(X,f)$ has $H$-shadowing if for all $\epsilon>0$ there exists $\delta>0$ such that if $(x_0,....,x_n)$ is a finite $\delta$-pseudo orbit, then there exists a $z \in X$ such that $d(f^i(z),x_i)<\epsilon$ and $f^n(z)=x_n$.
\end{definition}

\begin{definition}\label{def:H-top}$H$-shadowing- Topological definition:
A dynamical system  $(X,f)$ has  $H$-shadowing if for all $U \in FOC(X)$ , there exists a $V \in FOC(X)$ such that for all  finite $V$-pseudo orbits $(x_i)$,  we can find $z \in X$ such that $f^i(z),x_i \in U_i$ and $f^n(z)=x_n$.
\end{definition}

$H$-shadowing arises naturally when studying shadowing in shift spaces. In a shift space, if a $\delta$-pseudo orbit $(x_n)$ is shadowed by a point $z$, then there exists an $N \in \mathbb{N} $ such that for all $n\geq N$ we have that $f^n(z)=x_{n}$. If a  dynamical systems has $H$-shadowing, this can be used to characterize their $\omega$-limit sets \cite{barwell2013characterizations}.

\begin{theorem}
    Let $(X,f)$ be a dynamical system where $X$ is a compact metric space. The topological and metric definitions of  $H$-shadowing are equivalent.
\end{theorem}

\begin{proof}
    Suppose the topological definition of  $H$-shadowing holds in $(X,f)$. Let $\epsilon>0$. Since $X$ is compact we can choose an open cover $U$ of $X$ which consists of a finite number of $\frac{\epsilon}{4}$ balls. Choose $V \in FOC(X)$ which satisfies \Cref{def:H-top} with respect to $U$. Let $\delta>0$ denote the Lebesgue number of $V$. Let $x_0,x_1,\dots,x_n$ be a finite $\delta$-pseudo orbit. Then for all $i \in \{0,1,\dots,n\}$ we have $d(f(x_{i-1}),x_i)< \delta$. Hence there exists a $V_k \in V$ such that $f(x_{i-1}),x_{i} \in V_k$. Hence this $\delta$ pseudo orbit is also a $V$ pseudo orbit. Choose a $z$ which satisfies \Cref{def:H-top}. As we have topological $H$-shadowing, we know that for all $i \in \{0,1,\dots,n \}$, there exists a $U_j \in U$ such that $f^i(x)\text{ and }(x_i)$ both belong to $U_j$ with $f^n(z)=x_n$. Since $U$ consists of $\frac{\epsilon}{4}$-balls, we know that $d(f^i(z),x_i)<\frac{\epsilon}{2}$. Hence the metric definition of $H$-shadowing holds in $(X,f)$.

    Now we consider the case where the metric definition of $H$-shadowing holds. Consider an finite open cover $U$ of $X$. Let $\epsilon>0$ denote the Lebesgue number of $U$. Choose $\delta>0$ which satisfies \Cref{def:$H$-metric} with respect to $\epsilon$. As $X$ is compact we can choose a finite open cover $V$ of $X$ which consists of $\frac{\delta}{2}$-balls. Let  $x_0,x_1,\dots,x_n$ be a finite $V$-pseudo orbit. Hence for all $i\in \{0,1,\dots,n\}$, we have that $f(x_{i-1}),x_i \in V_k$ for some $V_k \in V$. As $V$ consists of $\frac{\delta}{2}$-balls, clearly $x_0,x_1,\dots,x_n$   is also a finite $\delta$-pseudo orbit. Choose $z$ which satisfies \Cref{def:$H$-metric}. $z$ satisfies $d(f^i(z),x_i)<\epsilon$ for $i \in \{0,1,\dots,n \}$ with $f^n(z)=x_n$. Hence there exists a $U_j$ such that $f^i(z),x_i \in U_j $. Hence $(X,f)$ also exhibits topological $H$-shadowing.

\end{proof}

\begin{definition}Almost lifting hit pseudo orbit property  (ALhP) - Topological definition:

We say a factor map  $\psi:X \longrightarrow Y $ has almost lifting hit pseudo orbit property if for all $V_{X}$ $\in$ $FOC(X)$, and $W_{Y}$ $\in$ $FOC(Y)$, there exists $V_{Y}$ $\in$ $FOC(Y)$,  such that if $y_0,....,y_n$ is a finite $V_{Y}$  pseudo orbit in Y, then there exists a  $V_{X}$ pseudo orbit in X, $x_0,....,x_n$ such that $\phi( x_{i}),y_{i}$ $\in$ $W_{k}$, for some $k$ and $\psi(x_n)=y_n$.

\end{definition}

\begin{definition}Almost lifting hit  pseudo orbit property (ALhP) - Metric definition:

We say a factor map  $\psi:X \longrightarrow Y $ has almost lifting hit pseudo orbit property if for all $\epsilon>0$ and $\eta>0$, there exists $\delta>0$  such that if $y_0,....,y_n$ is a finite $\delta$ pseudo orbit in Y, then there exists a  $\eta$ pseudo orbit in X, $x_0,....,x_n$ such that $d(\phi( x_{i}),y_{i})<\epsilon$, and $\psi(x_n)=y_n$.

\end{definition}

\begin{theorem}

Let $(X,f)$ and $(Y,g)$ be two dynamical systems. suppose $(X,f)$ has $H$-shadowing and $\psi:X \longrightarrow Y$ has the ALhP property, then $(Y,g)$ has $H$-shadowing.

\end{theorem}

\begin{proof}

Fix an open cover $U_{Y}$ of $Y$. Let $W_{Y}$ be an open cover of $Y$ such that if $W$,$W'$ $\in$  $W_{Y}$ and $W$ $\cap$ $W'$$\neq$$\emptyset$,  there exist $U$$\in$$U_{Y}$ such that  $W$ $\cup$ $W'$ $\subset U $. Let $U_{X}$ be an open cover of X such that $\psi(U_{X})$ refines $W_Y$. Since $(X,f)$ has h-shadowing , let $V_{X}$ witness this with respect to $U_{X}$. Since $\psi$ has ALhP, let $V_Y$ witness this with respect to $W_Y$ and $V_X$.

\bigskip

Let $y_0,....,y_n$  be a finite $V_Y$  pseudo orbit. Let $x_0,....,x_n$ be the  $V_X$ pseudo orbit, such that   $\psi((x_i))$ $W_Y$-shadows $(y_i)$ and $\psi(x_n)=y_n$.

\bigskip

Since $(X,f)$ has $H$-shadowing, we know that $x_i$ is $U_X$ shadowed by a point say z $\in$ X, such that $f^n(z)=x_n$. Hence for all $i$ we have that $\psi(f^i(z))$, $\psi(x_i)$ $\in$ $\psi(U_i)$ for some $U_i$ $\in$ $U_X$. Since $\psi$ is a factor map we know that $\psi(f^i(z))=g^i(\psi(z))$. And since $\psi(U_X)$ refines $W_Y$, it follows that there exists $W_i$ such that $g^i(\psi(z))$ and $\psi(x_i)$ lie in $W_i$. But we also know that $\psi(x_i)$ and $y_i$ belong to $W_{j}$ for some $j$. Hence $\psi(x_i)$ is in $W_{i}$ $\cap$ $W'_{j}$. But by construction of $W_Y$ we know there exists a $U_i$ such that $g^i(\psi(z))$ and $y_i$ lie in it. Hence $y_i$ is $U_Y$ shadowed by $\psi(z)$, and $g^n(\psi(z))=\psi(f^n(z))=\psi(x_n)=y_n$. Hence $(Y,g)$ has $H$-shadowing.

\end{proof}

\begin{theorem}

Suppose $\psi$ is a factor map from  $(X,f)$ to $(Y,g)$. We will prove that if $(Y,g)$ has $H$-shadowing, then $\psi$ must have the ALhP property. 

\end{theorem}

\begin{proof}

Let $\epsilon > 0$ and $\eta>0$. Choose $\delta >0$ such that $(Y,g)$  witnesses  $H$-shadowing with respect to $\epsilon$. Let $y_0,....,y_n$ be a finite $\delta$ pseudo orbit. Then it must be $\epsilon$-shadowed by a point say $z$, such that $g^n(z)=y_n$. Since $\psi$ is surjective we know there exists $x$ such that $\psi(x)=z$. Now consider $x,f(x),...,f^n(x)$. This is trivially a $\eta$ pseudo orbit. Note that for each $i$ we must have that $\psi(f^i(x))=g^i(\psi(x))=g^i(z)$, so  $\psi(f^i(x))$ $\epsilon$-shadows $y_0,....,y_n$. In particular $\psi(f^n(x))=g^n(\psi(x))=g^n(z)=y_n$. And so $(Y,g)$ has the ALhP property.

\end{proof}

In \cite{MT}, it was shown that $H$-shadowing is preserved under the induced hyperspace.

\section{Periodic shadowing}

\begin{definition}\label{def:2}Periodic shadowing - Metric definition:
    We say that $(X,f)$ has periodic shadowing if for all $\epsilon>0$ there exists $\delta>0$ such that if $(x_1,....,x_n)$ is a periodic $\delta$-pseudo orbit, then there exists a periodic $z \in X$ such that $d(f^i(z),x_i)<\epsilon$.
\end{definition}

\begin{definition}\label{def:1}  Periodic shadowing- Topological definition:
A dynamical system  $(X,f)$ has  periodic shadowing if for all $U \in FOC(X)$ , there exists a $V \in FOC(X)$ such that if $(x_i)$ is a periodic  $V$-pseudo orbits $(x_i)$,  we can find a periodic $z \in X$ such that $f^i(z),x_i \in U_i$.
\end{definition}

Periodic shadowing considers pseudo-orbits that approximate the behaviour of periodic orbits. It states that periodic pseudo-orbits can be shadowed by a genuine periodic point. It has been shown that there is a close relationship between a dynamical system exhibiting periodic shadowing, and it being $\Omega$-stable \cite{osipov2010periodic,chu2016periodic}.

\begin{theorem}
    Let $(X,f)$ be a dynamical system where $X$ is a compact metric space. The topological and metric definitions of periodic shadowing are equivalent.
\end{theorem}

\begin{proof}
    Suppose the topological definition of periodic shadowing holds in $(X,f)$. Let $\epsilon>0$. Since $X$ is compact we can choose an open cover $U$ of $X$ which consists of a finite number of $\frac{\epsilon}{4}$ balls. Choose $V \in FOC(X)$ which satisfies \Cref{def:1} with respect to $U$. Let $\delta>0$ denote the Lebesgue number of $V$. Let $x_0,x_1,\dots,x_n$ be a $\delta$-periodic pseudo orbit. Then for all $i \in \{0,1,\dots,n\}$ we have $d(f(x_{i-1}),x_i)< \delta$. Hence there exists a $V_k \in V$ such that $f(x_{i-1}),x_{i} \in V_k$.       Hence this $\delta$-periodic pseudo orbit is also a $V$-periodic pseudo orbit. Choose a periodic $z$ which satisfies \Cref{def:1}. As we have topological shadowing, we know that for all $i \in \{0,1,\dots,n \}$, there exists a $U_j \in U$ such that $f^i(x)\text{ and }(x_i)$ both belong to $U_j$. Since $U$ consists of $\frac{\epsilon}{4}$-balls, we know that $d(f^i(z),x_i)<\frac{\epsilon}{2}$. Hence the metric definition is satisfied.

    Now we consider the case where the metric definition holds. Consider an finite open cover $U$ of $X$. Let $\epsilon>0$ denote the Lebesgue number of $U$. Choose $\delta>0$ which satisfies \Cref{def:2} with respect to $\epsilon$. As $X$ is compact we can choose a finite open cover $V$ of $X$ which consists of $\frac{\delta}{2}$-balls. Let  $x_0,x_1,\dots,x_n$ be a $V$-periodic pseudo orbit. Hence for all $i\in \{0,1,\dots,n\}$, we have that $f(x_{i-1}),x_i \in V_k$ for some $V_k \in V$. As $V$ consists of $\frac{\delta}{2}$-balls, clearly $x_0,x_1,\dots,x_n$   is also a $\delta$-periodic pseudo orbit. Choose $z$ which satisfies \Cref{def:2}. $z$ is periodic and satisfies $d(f^i(z),x_i)<\epsilon$ for $i \in \{0,1,\dots,n \}$. Hence there exists a $U_j$ such that $f^i(z),x_i \in U_j $. Hence the topological definition is also satisfied.

\end{proof}

\begin{definition} Almost Lifting Periodic Pseudo orbit property (ALPP) Metric definition :
We say a factor map has ALPP if for all  $\epsilon > 0$ and for all $\eta > 0$, there exist a $\delta >0$, such that if $(y_i)_{i \in \omega}$ is a $\delta$ periodic pseudo orbit in Y, then there exists a $\eta $ periodic pseudo orbit in X, such that $d(\phi( x_{i}),y_{i} ) < \epsilon$.
\end{definition}

We can extend this definition more generally to topological spaces by using the following definition:

\begin{definition} Almost Lifting Periodic Pseudo orbit property (ALPP) Topological definition: We say a factor map has ALPP if for all  $V_{X}$ $\in$ $FOC(X)$, and $W_{Y}$ $\in$ $FOC(Y)$, there exists $V_{Y}$ $\in$ $FOC(Y)$,  such that if $(y_i)_{i \in \omega}$ is a $V_{Y}$ periodic pseudo orbit in Y, then there exists a  $V_{X}$  periodic pseudo orbit in X, such that $\phi( x_{i}),y_{i}$ $\in$ $W_{k}$, for some $k$. In this definition $FOC$ just means finite open cover.
\end{definition}

\bigskip

\begin{theorem}
Let $(X,f)$ and $(Y,g)$ be two dynamical systems. suppose $(X,f)$ has periodic shadowing and $\psi:X \longrightarrow Y$ has the ALPP property, then $(Y,g)$ has periodic shadowing.
\end{theorem}

\smallskip

\begin{proof} 
Fix an open cover $U_{Y}$ of $Y$. Let $W_{Y}$ be an open cover of $Y$ such that if $W$,$W'$ $\in$  $W_{Y}$ and $W$ $\cap$ $W'$$\neq$$\emptyset$,  there exist $U$$\in$$U_{Y}$ such that  $W$ $\cup$ $W'$ $\subset U $. Let $U_{X}$ be an open cover of X such that $\psi(U_{X})$ refines $W_Y$. Since (X,f) has periodic shadowing , let $V_{X}$ witness this with respect to $U_{X}$. Since $\psi$ has ALPP, let $V_Y$ witness this with respect to $W_Y$ and $V_X$.

\smallskip

Let  $(y_i)_{i \in \omega}$ be a periodic $V_Y$ periodic pseudo orbit. Let  $(x_i)_{i \in \omega}$ be the periodic $V_X$ pseudo orbit, such that   $\psi((x_i))_{i \in \omega}$ $W_Y$ shadows $(y_i)_{i \in \omega}$.

\smallskip

Since (X,f) has periodic shadowing, we know that $x_i$ is $U_X$ shadowed by a periodic point say z $\in$ X. Hence for all $i$ we have that $\psi(f^i(z))$, $\psi(x_i)$ $\in$ $\psi(U_i)$ for some $U_i$ $\in$ $U_X$. Since $\psi$ is a factor map we know that $\psi(f^i(z))=g^i(\psi(z))$. And since $\psi(U_X)$ refines $W_Y$, it follows that there exists $W_i$ such that $g^i(\psi(z))$ and $\psi(x_i)$ lie in $W_i$. But we also know that $\psi(x_i)$ and $y_i$ belong to $W_{j}$ for some $j$. Hence $\psi(x_i)$ is in $W_{i}$ $\cap$ $W'_{j}$. But by construction of $W_Y$ we know there exists a $U_i$ such that $g^i(\psi(z))$ and $y_i$ lie in it. Hence $y_i$ is $U_Y$ shadowed by $\psi(z)$, and since the image of periodic points under factor maps are periodic, we have shown that $(Y,g)$ has periodic shadowing>

\end{proof}

\begin{theorem}

Suppose $\psi$ is a factor map from  $(X,f)$ to $(Y,g)$. We will prove that if $(Y,g)$ has the periodic shadowing property, then $\psi$ must have the ALPP property. 
\end{theorem}

\begin{proof} 

Let $\epsilon> 0$ and $\eta > 0$. Choose $\delta > 0$ such that it witnesses periodic shadowing with respect to $\epsilon$. Let $(y_i)$ be a $\delta$ periodic pseudo orbit. Then it must be shadowed by a periodic point say $z$. So there exists a $n$ such that $g^n(z) = z$. Consider a covering of X by open balls of radius $\frac{\eta}{2}$. Since X is compact this covering has a finite sub cover.
Hence $\psi(z)^{-1}$ is contained in finite number of these open balls of radius $\frac{\eta}{2}$. Since $\psi$ is surjective we know there exists $x^*$
such that $\psi(x^*)$ = z. Consider the sequence $x^*, f(x^*), f^2(x^*), \dots$ Note that $x^*, f^n(x^*) f^{2n}(x^*), f^{3n}(x^*), \dots$ all belong in $\psi^{-1}(z)$. By the
pigeonhole principle, at least one of the open balls of radius $\frac{\eta}{2}$ that cover $\psi^{-1}(z)$ contains $f^{in}(x^*),f^{kn}(x^*)$  where $i < k$. But then consider the $\eta$ periodic pseudo orbit given by $f^{in}(x),f^{in+1}(x),\dots f^(kn-1)(x^*), f^{in}(x),f^{in+1}(x),\dots f^(kn-1)(x^*), f^{in}(x^*),f^{in+1}(x^*),\dots$

Clearly this $\eta$ periodic pseudo orbit satisfies the condition needed for the ALPP property to hold. Hence the map $\psi$ satisfies the ALPP property.

\end{proof}

\section{Ramsey shadowing}

\begin{definition} Furstenberg family:

We call $\mathcal{F}$ a Furstenberg family, if it is a  nonempty set of subsets of $\mathbb{N}$ and it is closed under supersets.

\end{definition}

\begin{definition}Ramsey property:

We say $\mathcal{F}$,a Furstenberg family, has the Ramsey property, if whenever $A \in \mathcal{F}$ and $A=\bigcup_{i<n} A_i$, then $A_i \in  \mathcal{F}$ for some i.

\end{definition}

Ramsey shadowing has shown itself to be strongly connected to minimal points and other fundamental concepts in dynamical systems, including recurrence. In particular, it has been shown that Ramsey shadowing is equivalent to the existence of a dense set of minimal points \cite{brian2016Ramsey}.

\begin{definition}\label{def:Ramsey-metric}Ramsey shadowing- Metric definition:
A dynamical system  $(X,f)$ has Ramsey shadowing if for all $\epsilon>0$, for all $\mathcal{F}$ with the Ramsey property, any sequence $(x_i)$ can be $\epsilon$-shadowed on a set in $\mathcal{F}$. 
\end{definition}

\begin{definition}\label{def:Ramsey-Top}Ramsey shadowing- Topological definition:
A dynamical system  $(X,f)$ has Ramsey shadowing if for all $U$ finite open cover of $X$, for all $\mathcal{F}$ with the Ramsey property, any sequence $(x_i)$ can be $U$-shadowed on a set in $\mathcal{F}$. 
\end{definition}

\begin{theorem}
    Let $(X,f)$ be a dynamical system where $X$ is a compact metric space. The topological and metric definitions of  Ramsey shadowing are equivalent.
\end{theorem}

\begin{proof}
    Suppose the topological definition of  Ramsey shadowing holds in $(X,f)$. Let $\epsilon>0$. Since $X$ is compact we can choose an open cover $U$ of $X$ which consists of a finite number of $\frac{\epsilon}{4}$ balls. Pick an arbitrary Furstenberg family with the Ramsey property and a arbitrary sequence $(x_i)_{i\in \mathbb{N}_0}$. Choose a $z$ which satisfies \Cref{def:Ramsey-Top}. As we have topological Ramsey shadowing, we know that there exists a $z \in X$ and a $A \in \mathcal{F}$ such that for all $i \in A$, there exists a $U_j \in U$ such that $f^i(z)\text{ and }(x_i)$ both belong to $U_j$. Since $U$ consists of $\frac{\epsilon}{4}$-balls, we know that $d(f^i(z),x_i)<\frac{\epsilon}{2}$. Hence the metric definition of Ramsey shadowing holds in $(X,f)$, as  $(x_i)_{i\in \mathbb{N}_0}$ is $\epsilon$-shadowed by $z$ on $A$.

    Now we consider the case where the metric definition of Ramsey-shadowing holds. Consider a finite open cover $U$ of $X$. Let $\epsilon>0$ denote the Lebesgue number of $U$.  Pick an arbitrary sequence  $(x_i)_{i\in \mathbb{N}_0}$  and a arbitrary Furstenberg family with the Ramsey property. Choose $z$ which satisfies \Cref{def:Ramsey-metric}. We know that there exists a $z \in X$ and a $B \in \mathcal{F}$ such that for all $i \in A$, we have $d(f^i(z),x_i)< \epsilon$. Therefore  $(x_i)_{i\in \mathbb{N}_0}$ is $U$-shadowed on $B$, and so $(X,f)$ also exhibits topological Ramsey shadowing.

\end{proof}

\begin{theorem}

Suppose $\psi$ is a factor map from  $(X,f)$ to $(Y,g)$.  We will show that if $(X,f)$ has Ramsey shadowing, then so does $(Y,g)$
\end{theorem}

\begin{proof}

Fix an open cover $U_{Y}$ of $Y$. Let $W_{Y}$ be an open cover of $Y$ such that if $W$,$W'$ $\in$  $W_{Y}$ and $W$ $\cap$ $W'$$\neq$$\emptyset$,  there exist $U$$\in$$U_{Y}$ such that  $W$ $\cup$ $W'$ $\subset U $. Let $U_{X}$ be an open cover of X such that $\psi(U_{X})$ refines $W_Y$. 

\bigskip

Let $(y_i) \in Y$. Let $\mathcal{F}$ be arbitrary with the Ramsey property. Let $(x_i)$ be the sequence defined by $x_i=\psi^{-1}(y_i)$. By Ramsey shadowing we know that $x_i$, is $U_X$-shadowed by a point $z \in X$ on some $A \in \mathcal{F}$. we note that $\psi(x_i)$ certainly $W_Y$-shadows $(y_i)$ on $A$. from now on we assume the subscript is on A. But then  $\psi(x_i)$  and $(y_i)$ belong to some $W_k$. But also $\psi(x_i)$ and $\psi(f^i(z))$ also belong to some $W_t$, hence $W_t$ and $W_k$ is non empty so by construction we must have that $\psi(f^i(z))$ and $(y_i)$ belong to some $U_s$. But note that  $\psi(f^i(z))=g^(\psi(z))$. And so $\psi(z)$ $U_Y$-shadows $(y_i)$ on A.

\end{proof}

\begin{theorem}

If the hyperspace $(2^X,2^f)$ has Ramsey shadowing then so does the $(X,f)$.

\end{theorem}

\begin{proof}
    Let $\epsilon>0$. Let $\mathcal{F}$ be an arbitrary Furstenberg family with the Ramsey property. Consider an arbitrary sequence $(x_i)$. Now consider the sets $\{x_1\},\{x_2\},...$. These are compact sets hence belong to $2^X$. Since the hyperspace has Ramsey shadowing we know that there exist a compact set $A \in 2^X$ which $\epsilon$-shadows ${x_1},{x_2},..$ on a set $B$ in $\mathcal{F}$. But then by the definition of Hausdorff distance we see that any element of $A$ must shadow $(x_i)$ on $B \in \mathcal{F}$. Hence $(X,f)$ has Ramsey shadowing.

\end{proof}

\begin{theorem}
The inverse limit of a sequence of dynamical systems, each exhibiting the Ramsey shadowing property and satisfying the Mittag-Leffler condition, also has the Ramsey shadowing property.

\end{theorem}

\begin{proof}

Suppose $(X,f)$ is conjugate to the inverse limit system $(X_n,f_n,g^{\gamma}_k)$ where each $(X_n,f_n)$ has the Ramsey shadowing property. Let $W$ be a finite open cover of $X$. Let $F$ be a Furstenberg family with the Ramsey property. We can find a $k$ and a finite open cover $B_k$ of $X_k$ such that the  $B^*=\{\pi^{-1}(B)\cap X: B\in B_k\}$ refine $W$. 
\newline

Using the Mittag-Leffler condition we know there exists a $\gamma$ which satisfies the condition with respect to $k$. $B_{\gamma}=\{{(g^{\gamma}_{k})}^{-1}(B): B\in B_k\}$. Let $(x_i)$ be an arbitrary sequence of point in $X$. This means that $(\pi_\gamma(x_i))$ is a sequence of points in $X_\gamma$. Hence we know that there exists an point $x_\gamma$ which $B_\gamma$-shadows on it on some $A \in \mathcal{F}$. Hence it follows that $g^\gamma_{k}(x^*)$ $B_k$-shadows $(\pi_{k}(x_i))$ on $A \in \mathcal{F}$. 
\newline 

Note that we have that $g^\gamma_{k}(x^*) \in g^\gamma_{k}({\ \bigcap(f^{-i}(B^{i}_{\gamma})) }\ ) \subseteq \bigcap\limits_{i \in A}(f^{-i}(B^{i}_{k})  $

And so, using the Mittag-Leffler condition we know there exists an $z$ such that $(z)_k \in g^\gamma_{k}({\ \cap(f^{-i}(B^{i}_{\gamma})) }\ ) $. And so we must have that $f^i(z),x_i \in \pi^{-1}_k(W^{i}_k) \cap X $. Hence $(x_i)$ is $B*$-shadowed on $A \in F$. But since $B^*$ refines $W$, we must also have that $(x_i)$ is $W$-shadowed on $A \in F$.

\end{proof}

\section{Ergodic Shadowing}

\begin{definition} Ergodic pseudo-orbit - Metric Definition:

Let $\delta>0$. Let $(x_i)$ be a sequence. Let $N((x_i),\delta)= \{i: d(f(x_i),x_{i+1}) \geq \delta    \}$. Let $N((x_i),\delta,n)= \{i: d(f(x_i),x_{i+1}) \geq \delta,  i=0,...,n    \}$. Let $M((x_i),\delta,z)= \{i: d(f^i(z),x_{i}) \geq \delta \}$ and $M((x_i),\delta,z,n)= \{i: d(f^i(z),x_{i}) \geq \delta,  i=0,...,n    \}$. We call $(x_i)$ a $\delta$-ergodic pseudo orbit provided that $\lim_{n\to\infty} \frac{ \lvert N((x_i),\delta,n) \rvert }{n}=0.$ We say that a $\delta$-ergodic pseudo orbit $(x_i)$ is $\epsilon$-ergodic shadowed by a point $z$ provided that $\lim_{n\to\infty} \frac{  \lvert M((x_i),\epsilon,z,n) \rvert }{n}=0.$ 

\end{definition}

\begin{definition}
     Ergodic pseudo-orbit - Topological definition:

Let $V\in FOC(X)$, we call $(x_i)$ a $V$-ergodic pseudo-orbit if:

$\lim_{n\to\infty} \frac{ \lvert \{i\in [0,...n]: \forall V_k \in V, \text{ } f(x_i) \notin V_k \text{ or } x_{i+1} \notin V_k  \} \rvert }{n}=0$

\end{definition}

\begin{definition} Ergodic shadowing - Metric definition:

We say the dynamical system has ergodic-shadowing if for all $\epsilon>0$, there exists $\delta>0$ such that if $(x_i)$ is a $\delta$-ergodic pseudo orbit, then it is $\epsilon$ ergodic shadowed by a point $z$.

\end{definition}

Ergodic shadowing was introduced in \cite{fakhari2010shadowing}, and has shown itself to be closely linked to  various notions of mixing. In fact, in \cite{fakhari2010shadowing}, they showed that ergodic shadowing was equivalent to having both the normal shadowing property as well as topological mixing.

\begin{definition}Ergodic shadowing- Topological definition:
A dynamical system  $(X,f)$ has ergodic shadowing if for all $U \in FOC(X)$ , there exists a $V \in FOC(X)$ such that for all ergodic $V$-pseudo orbits, we can find a $z \in X$ such that : 

$\lim_{n\to\infty} \frac{ \lvert \{i\in [0,...n]: \forall U_k \in U, \text{ } f^i(z) \notin U_k \text{ or } x_i \notin U_k  \} \rvert }{n}=0$
\end{definition}

As in the previous cases, the metric and topological definitions of the Ergodic shadowing can be shown to be equivalent. We omit the proof, as it follows essentially the same arguments as before.

\begin{definition}Ergodic Lifting Property (ELP) - Topological definition:

 We say a factor map has ELP if for all $U_X \in FOC(X)$ and for all $W_Y \in FOC(Y)$, there exists $V_Y \in FOC(Y)$, such that if $(y_i)$ is a  $V_Y$-ergodic pseudo orbit, there exists a $U_X$-ergodic pseudo orbit $(x_i)$ in X such that :

 $\lim_{n\to\infty} \frac{ \lvert \{i\in [0,...n]: \forall W_{Y_t} \in W_Y \psi(x_i) \notin W_{Y_t} \text{ or } y_i \notin W_{Y_t} \} \rvert }{n}=0$.

\end{definition}

\begin{definition}Ergodic Lifting Property (ELP) - Metric definition:

 We say the dynamical system has ELP if for all $\epsilon>0$ and for all $\eta>0$, there exists $\delta>0$, such that if $(y_i)$ $\delta$-ergodic pseudo orbit, there exists a $\eta$-ergodic pseudo orbit $(x_i)$ in X such that :

 $\lim_{n\to\infty} \frac{\lvert \{i\in [0,...n]: d(\psi(x_i),y_i)>\epsilon\} \rvert }{n}=0$.

\end{definition}

\begin{theorem}
    Suppose $\psi$ is a factor map from  $(X,f)$ to $(Y,g)$. We will prove that if $(Y,g)$ has ergodic shadowing, then $\psi$ must have ELP.

\end{theorem}

\begin{proof}

Let $\epsilon > 0$  and $\eta>0$. Choose $\delta>0$, such that $Y$ satisfies ergodic shadowing with respect to $\epsilon$. Let $(y_i)$ be a $\delta$-ergodic pseudo orbit. We know it is $\epsilon$-ergodic shadowed by a  point say $z$. 
\newline

So we must have that $\lim_{n\to\infty} \frac{ \lvert  M((x_i),\epsilon,z,n) \rvert }{n}=0$. Take $x \in \psi^{-1}(z)$ and consider the sequence $(x_i)$ given by its orbit $(x,f(x),...)$. This is clearly a $\eta$-ergodic pseudo orbit. And by semi-conjugacy we know that when $\psi(f^n(x))=g^n(z)$. And since $z$ $\epsilon$-ergodic shadows $y_i$ we know that:

 $\lim_{n\to\infty} \frac{ \lvert \{i\in [0,...n]: d(\psi(x_i),y_i)>\epsilon\} \rvert }{n}=0$. Hence $\psi$ must have ELP.

\end{proof}

\begin{theorem}

Let $(X,f)$ and $(Y,g)$ be two dynamical systems. suppose $(X,f)$ has ergodic shadowing and $\psi:X \longrightarrow Y$ has ELP , then $(Y,g)$ has ergodic shadowing.
\end{theorem}

\begin{proof}

Fix an open cover $U_{Y}$ of $Y$. Let $W_{Y}$ be an open cover of $Y$ such that if $W$,$W'$ $\in$  $W_{Y}$ and $W$ $\cap$ $W'$$\neq$$\emptyset$,  there exist $U$$\in$$U_{Y}$ such that  $W$ $\cup$ $W'$ $\subset U $. Let $U_{X}$ be an open cover of X such that $\psi(U_{X})$ refines $W_Y$. Since $X$ has ergodic-shadowing, let $V_X$ be an open cover that satisfies this with respect to $U_X$. Since $\psi$ has ELP let $V_Y$ respect it with regards to $W_Y$ and $V_X$.

\bigskip

Let $(y_i) \in Y$ be a $V_Y$-ergodic pseudo orbit. Let $(x_i)$ be the $V_X$-ergodic pseudo orbit in $X$ that satisfies ELP. By ergodic-shadowing we know that $x_i$, is $U_X$-shadowed by a point $z \in X$, in such a way that :

\bigskip

$\lim_{n\to\infty} \frac{ \lvert \{i\in [0,...n]: \forall A \in U_{X}, \textbf{ }f^i(z)\notin A \text{ or } x_i \notin A  \} \rvert }{n}=0$

\bigskip

 But by the ELP property we note that $\psi(x_i)$ $W_Y$-shadows $(y_i)$ in such a way that:

 \bigskip
 
 $\lim_{n\to\infty} \frac{ \lvert \{i\in [0,...n]: \forall B \in W_{Y}, \text{ }  y_i \notin B \text{ or } \psi(x_i) \notin B  \} \rvert }{n}=0$

 \bigskip

Since $\psi$ refines $W_Y$, we know that whenever  $(x_i)$ and $(f^i(z))$  belong the same $U_{X_k}$, then $\psi(x_i)$ and $\psi(f^i(z))=g^i(\psi(z))$ must also belong to some $W_t$.

We want to show that :

 \bigskip
 
 $\lim_{n\to\infty} \frac{ \lvert \{i\in [0,...n]: \forall B \in W_{Y},\text{ } y_i \notin B \text{ or } \psi(f^i(z)) \notin B  \} \rvert }{n}=0$.

 Let $\epsilon>0$. We know there exists $M \in N$ such that for all $n>M$:

 \bigskip

$ \frac{ \lvert \{i\in [0,...n]: \forall A\in U_{X}, f^i(z) \notin A \text{ or } x_i \notin A \} \rvert }{n}< \frac{\epsilon}{2}$ 

\bigskip
and
\newline

$ \frac{ \lvert  \{i\in [0,...n]: \forall B \in W_{Y}, \text{ } y_i \notin B \text{ or } \psi(x_i) \notin B \} \rvert }{n}< \frac{\epsilon}{2}$.
\bigskip

We note that $\{i\in [0,...n]: \forall B \in U_{Y}, \text{ } y_i \notin B \text{ or } \psi(f^i(z)) \notin B \} \subseteq \{i\in [0,...n]: \forall A \in U_{X}, \text{ } x_i \notin A \text{ or } f^i(z) \notin A \}  \cup    \{i\in [0,...n]: \forall B \in W_{Y}, \text{ } y_i \notin B \text{ or } \psi(x_i) \notin B $ 
Hence we must have that:

\bigskip

 $ \frac{ \lvert \{i\in [0,...n]: \forall B \in U_{Y}, \text{ } y_i \notin B \text{ or } \psi(f^i(z)) \notin B  \} \rvert }{n}< \frac{\epsilon}{2} + \frac{\epsilon}{2} = \epsilon   $. 
\bigskip

 Hence:
 $\lim_{n\to\infty} \frac{ \lvert \{i\in [0,...n]: \forall B \in U_{Y}, \text{  } y_i \notin B \text{ or } \psi(f^i(z)) \notin B \} \rvert }{n}=0$.
\bigskip

 But then $\psi(f^i(z))=g^i(\psi(z))$. And so $\psi(z)$ $U_Y$-ergodic shadows $(y_i)$, hence we have ergodic shadowing.

\end{proof}

\begin{theorem}

If the hyperspace $(2^X,2^f)$ has  ergodic shadowing then so does the $(X,f)$.

\end{theorem}

\begin{proof}
        Let $\epsilon>0$. Choose $\delta$ such that $(2^X,2^f)$ has ergodic shadowing. Now consider a $\delta$-ergodic pseudo orbit in X, $(x_i)$. Now consider the sets $\{x_1\},\{x_2\},...$. These are compact sets hence belong to $2^X$. Since the hyperspace has ergodic-shadowing we know that there exist a compact set $A \in 2^X$ such that :
    \newline
    \newline
    $\lim_{n\to\infty} \frac{   M(\{x_i\},\epsilon,(2^f)^i(A),n) }{n}=0.$     \newline

    But then by the definition of Hausdorff distance we see that any element of $a \in A$ must satisfy:

        $\lim_{n\to\infty} \frac{   M(x_i,\epsilon,f^i(a),n) }{n}=0.$
    Hence $(X,f)$ also has ergodic shadowing.
        
\end{proof}

\section{Density shadowing}

\begin{definition}\label{def:density-metric}  \underline{d}-shadowing- Metric definition :

We say the dynamical system has \underline{d}-shadowing if for all $\epsilon>0$, there exists $\delta>0$ such that if $(x_i)$ is a $\delta$-ergodic pseudo orbit, then there is a point $z$ such that :

$\liminf_{n\to\infty} \frac{ \lvert \{i\in [0,...n]: d(f^i(z),x_i)<\epsilon\} \rvert }{n}>0$

\end{definition}

Density shadowing was introduced  in \cite{fakhari2010subshadowings} and has been shown to be a useful tool in the study of equicontinuity and transitivity \cite{Wu2014OnVD}.

\begin{definition}\label{def:density-top} \underline{d}-shadowing- Topological definition:
A dynamical system  $(X,f)$ has \underline{d} shadowing if for all $U \in FOC(X)$ , there exists a $V \in FOC(X)$ such that for all ergodic $V$-pseudo orbits, we can find a $z \in X$ such that : 

$\liminf_{n\to\infty} \frac{ \lvert \{i\in [0,...n]: f^i(z) \text{ and } x_i \in  U_k \text{ for some k }\} \rvert }{n}>0$

\end{definition}

\begin{theorem}
    Let $(X,f)$ be a dynamical system where $X$ is a compact metric space. The topological and metric definitions of  density shadowing are equivalent.
\end{theorem}

\begin{proof}
    Suppose the topological definition of  Density-shadowing holds in $(X,f)$. Let $\epsilon>0$. Since $X$ is compact we can choose an open cover $U$ of $X$ which consists of a finite number of $\frac{\epsilon}{4}$ balls. Choose $V \in FOC(X)$ which satisfies \Cref{def:density-top} with respect to $U$. Let $\delta>0$ denote the Lebesgue number of $V$.  $N((x_i),\delta,n)= \{i: d(f(x_i),x_{i+1}>\delta,  i=0,...,n    \}$.  $N((x_i),V,n)= \{i:\nexists\text{ } $j$ \text{ such that } f(x_i),x_{i-1} \in V_j,  i=0,...,n    \}$.\

    Let $(x_i)_{i\in\mathbb{N}_0}$ be a $\delta$-ergodic pseudo orbit. Then  $\lim_{n\to\infty} \frac{ \lvert N((x_i),\delta,n) \rvert }{n}=0$. Pick $k \in N((x_i),V,n) $ , then there is no $V_j \in V $ such that   $f(x_k),x_{k+1} \in V_j$. As $\delta$ is the Lebesgue number of $V$ this implies that $d(f(x_k),x_{k+1})> \delta$ and so $ N((x_i),V,n) \subseteq N((x_i),\delta,n)$. Hence $\lim_{n\to\infty} \frac{ \lvert N((x_i),V,n) \rvert }{n}=0$ and so this $\delta$ ergodic pseudo orbit is also a $V$-ergodic pseudo orbit. Choose a $z$ which satisfies \Cref{def:density-top}. As we have topological density-shadowing, we know that $\liminf_{n\to\infty} \frac{ \lvert \{i\in [0,...n]: d(f^i(z),x_i) \in U \} \rvert }{n}>0$. Since $U$ consists of $\frac{\epsilon}{4}$-balls, we know that whenever $f^i(z),x_i \in U_j$ we must also have that $d(f^i(z),x_i)<\frac{\epsilon}{2}$. Hence $\liminf_{n\to\infty} \frac{ \lvert \{i\in [0,...n]: d(f^i(z),x_i)< \epsilon \} \rvert }{n}>0$, and the metric definition of density shadowing also holds.

    Now we consider the case where the metric definition of density shadowing holds. Consider an finite open cover $U$ of $X$. Let $\epsilon>0$ denote the Lebesgue number of $U$. Choose $\delta>0$ which satisfies \Cref{def:density-metric} with respect to $\epsilon$. As $X$ is compact we can choose a finite open cover $V$ of $X$ which consists of $\frac{\delta}{2}$-balls. Let  $(x_i)_{i \in \mathbb{N}_0}$ be a $V$ ergodic pseudo orbit. Then $\lim_{n\to\infty} \frac{ \lvert N((x_i),V,n) \rvert }{n}=0$.    As $V$ consists of $\frac{\delta}{2}$-balls, clearly $\{i: f(x_i),x_{i+1} \in V,  i=0,...,n    \}   \subseteq \{i: d(f(x_i),x_{i+1} < V,  i=0,...,n    \}$. By taking complements we see that $N((x_i),\delta,n) \subseteq N((x_i),V,n)$, and so $\lim_{n\to\infty} \frac{ \lvert N((x_i),\delta,n) \rvert }{n}=0$. Hence $(x_i)_{i \in \mathbb{N}_0}$ is also a $\delta$-ergodic pseudo orbit. Choose $z$ which satisfies \Cref{def:density-metric}. Then $\liminf_{n\to\infty} \frac{ \lvert  \{i\in [0,...n]: d(f^i(z),x_i)< \epsilon \} \rvert }{n}>0$. As $\epsilon$ is the Lebesgue number of $U$, we know that whenever $d(f^i(z),x_i)< \epsilon$ we must also have that $f^i(z),x_i \in U_j$ for some $U_j \in U$. This implies that $\liminf_{n\to\infty} \frac{ \lvert \{i\in [0,...n]: f^i(z),x_i \in U\} \rvert }{n}>0$ and so $(X,f)$ also exhibits topological $\underline{d}$-shadowing.

\end{proof}

\begin{theorem}
    If the hyperspace $(2^X,2^f)$ has $\underline{d}$-shadowing then so does the $(X,f)$.

\end{theorem}
\begin{proof}

    Let $\epsilon>0$. Choose $\delta$ such that $(2^X,2^f)$ has $\underline{d}$-shadowing. Consider an ergodic $\delta$ pseudo orbit $(x_i)$. Now consider the sets $\{x_1\},\{x_2\},...$. These are compact sets hence belong to $2^X$. Since the hyperspace has $\underline{d}$-shadowing we know that there exist a compact set $A \in 2^X$ such that :
    \newline
    \newline
    $\liminf_{n\to\infty} \frac{ \lvert \{i\in [0,...n]: d((2^f)^i(A),\{x_i\})<\epsilon\} \rvert }{n}>0$
    \newline

    But then by the definition of Hausdorff distance we see that any element of $a \in A$ must satisfy:

        $\liminf_{n\to\infty} \frac{ \lvert \{i\in [0,...n]: d(f^(a),\{x_i\})<\epsilon\} \rvert }{n}>0$.

    Hence $(X,f)$ has $\underline{d}$-shadowing.

\end{proof}

\begin{theorem}
The inverse limit of a sequence of dynamical systems, each exhibiting  $\underline{d}$-shadowing and also satisfying the Mittag-Leffler condition, also has  $\underline{d}$ shadowing.
\end{theorem}

\begin{proof}

Suppose $(X,f)$ is conjugate to the inverse limit system $(X_n,f_n,g^{\gamma}_k)$ where each $(X_n,f_n)$ has the $\underline{d}$ shadowing property. Let $W$ be a finite open cover of $X$. We can find a $k$ and a finite open cover $B_k$ of $X_k$ such that the  $B^*=\{\pi^{-1}(B)\cap X: B\in B_k\}$ refine $W$. 
\newline

Using the Mittag-Leffler condition we know there exists a $\gamma$ which satisfies the condition with respect to $k$. Consider the open cover of $ W_ \gamma$ given by $B_{\gamma}=\{{(g^{\gamma}_{k})}^{-1}(B): B\in B_k\}$. Since $(X_\gamma, f_\gamma)$ has \underline{d}- shadowing, let $V_\gamma \in FOC(X_\gamma)$ be the open cover that satisfies the definition. Let $V= \{ \pi^{-1}_{\gamma}(V)\cap X: V\in V_\gamma \} $. Let $(x_i)$ be an ergodic $V$-pseudo orbit sequence in $X$, which is given by the sequence $(  \pi^{-1}(V^i)\cap X  )$. This means that $(\pi_\gamma(x_i))$ is an ergodic $V_\gamma$-pseudo orbit in $X_\gamma$ with sequence $(V^{i}_\gamma)$.  Since $(X_\gamma, f_\gamma)$ has \underline{d}- shadowing  we know that there exists an point $x_\gamma$ such that:
\newline
\newline
$\liminf_{n\to\infty} \frac{ \lvert \{i\in [0,...n]: d(f^i(z),x_i) \in B_{\gamma} \} \rvert }{n}>0$
\newline

Hence it follows that $g^\gamma_{k}(x^*)$ satisfies :
\newline
\newline
$\liminf_{n\to\infty} \frac{ \lvert \{i\in [0,...n]: d( g^{\gamma}_{k}(z) , (x_i)_k) \in B_k \} \rvert }{n}>0$
\newline

Let $A^*$ denotes the values for which $d(f^i(z),x_i) \in B_{\gamma}$.
Note that we have that $g^\gamma_{k}(x^*) \in g^\gamma_{k}({\ \bigcap\limits_{i \in A^*} (f^{-i}(B^{i}_{\gamma})) }\ ) \subseteq \cap(f^{-i}(B^{i}_{k})  $.

And so, using the Mittag-Leffler condition we know there exists an $z$ such that $(z)_k \in g^\gamma_{k}({\ \bigcap\limits_{i \in A^*}(f^{-i}(B^{i}_{\gamma})) }\ ) $. And so we must have that $f^i(z),x_i \in \pi^{-1}_k(W^{i}_k) \cap X $ for $i \in A^*$. 

Hence $(x_i)$ is $\underline{d}$-shadowed. But since $B^*$ refines $W$, we must also have that $(x_i)$ is $\underline{d}$-shadowed.

\end{proof}

\section*{Acknowledgments}

The first author, Bura, is thankful to EPSRC grants EP/S02297X/1 for financial support during the development of this article.

\bibliographystyle{alpha}
\bibliography{ref1}

\end{document}